\documentclass[10pt,a4paper]{article}
\usepackage[latin1]{inputenc} 
\usepackage[main=english,spanish]{babel}
\usepackage[T1]{fontenc} 
\usepackage[all]{xy}
\usepackage{fancyhdr}
\usepackage[dvips]{geometry}
\usepackage{indentfirst,amsfonts,amsmath,amsbsy,amsthm,amssymb,amscd}
\usepackage{enumerate}
\usepackage[normalem]{ulem}  
\usepackage{makeidx} 
\usepackage{color}
\usepackage{ mathrsfs }
\usepackage{dsfont}
\usepackage{hyperref}
\usepackage{makeidx}
\usepackage{calligra}
\usepackage[affil-it]{authblk}
\usepackage{graphicx}

\catcode`@=11 \@addtoreset{equation}{section} \catcode`@=12

\def\XXint#1#2#3{{\setbox0=\hbox{$#1{#2#3}{\int}$ }
\vcenter{\hbox{$#2#3$ }}\kern-.6\wd0}}

\newtheorem{theorem}{Theorem}[section]

\numberwithin{equation}{section}

\newcommand{\Beq}{\begin{equation}}
\newcommand{\Eeq}{\end{equation}}
\newcommand{\BS}{\begin{subequations}}
\newcommand{\ES}{\end{subequations}}
\newcommand{\Beqn}{\begin{equation*}}
\newcommand{\Eeqn}{\end{equation*}}
\newcommand{\Beqa}{\begin{eqnarray}}
\newcommand{\Eeqa}{\end{eqnarray}}
\newcommand{\Beqan}{\begin{eqnarray*}}
\newcommand{\Eeqan}{\end{eqnarray*}}
\newcommand{\N}{\hbox{$ I\kern -0.23em N$}}

\newcommand{\R}{{\mathbb R}}

\newcommand{\eps}{\varepsilon}

\newcommand\supp{\mathop{\rm supp}}

\begin{document}
\title
{Equivalent formulations for steady periodic water waves of fixed mean-depth 
 with discontinuous vorticity}

\author{Silvia Sastre-Gomez}

\affil{ Department of Applied Mathematics, School of Mathematical 
Sciences, 
\\University College Cork, Ireland.}

\maketitle

\begin{abstract}
In this work we prove the equivalence between three different weak formulations of the
steady periodic water wave problem where the vorticity is discontinuous. In particular, we prove 
that generalised versions of the standard Euler and stream function formulation of the governing 
equations are equivalent to a weak version of the recently introduced modified-height formulation. 
The weak solutions of these formulations are considered in H\"older 
spaces. \\
\\
\end{abstract}

\section{Introduction}
In this paper we consider steady periodic water waves, which propagate over a flat bed with a 
specified 
mean-depth, and which have discontinuous vorticity distribution. In particular, we prove the 
equivalence between three different weak formulations of the governing equations, namely the 
generalised standard Euler equation and stream function formulations, and the modified-height 
formulation. The standard governing equations for perfect (inviscid and incompressible) fluids are 
given by the Euler equation together with associated boundary conditions. Often it proves useful 
to reformulate these equations in terms of a stream function, leading to a semilinear elliptic 
equation with nonlinear boundary conditions. Both of these formulations 
are free boundary problems, and an inherent difficulty in their solution is the determination of the 
wave's free surface,  cf. \cite{Cons_book, ConJPA,J90,Toland_96}. 

One way to by-pass this difficulty is to employ a semi-hodograph change of variables to transform 
to a fixed-domain, with the trade-off being that our PDE system becomes quite more involved and 
complicated than previously. The standard transformation which is typically employed to this end 
in studying waves with vorticity (which model wave-current interactions \cite{Cons_book, 
ConJPA,TK97}) is the Dubreil-Jacotin transformation  \cite{DubJac}, whereby the system of 
governing equations may then be expressed in terms of a height function. This approach has 
been successfully implemented in \cite{CS2,CS11} in using local and global bifurcation theory to 
prove the existence of steady rotational water waves with small and large amplitude. Motivated by 
this work, there have been an extensive analytical studies of periodic waves with vorticity 
\cite{Con_Esc_deep_04, Cons_Esch_04, Cons_Esch_11, CS2, Con_Var_11, EMM1, 
Hen_anc_14, Hen_Bog_13,Hen_Bog_15, Var_09, Var_11,Wa6b,Wah3,Wa1}.

In this paper,  we are interested in a recently developed modified-height function reformulation of 
the governing equations, which follows when a variation on the Dubreil-Jacotin transformation is 
invoked. The formulation was employed in \cite{HDisp,HDisp13, Hen13, DH1} in order to prove 
the existence of rotational water waves of a fixed mean-depth $d$, as opposed to the approach in 
\cite{CS11} where the mass flux 
$p_0$ instead is fixed. Fixing the mean-depth of the wave is heuristically and physically quite 
natural, since the 
height can be measured more easily than the mass flux. Moreover, in \cite
{KoS1} it was noticed that fixing the mass flux does not fix the depth,  since it was observed 
numerically that on a bifurcation curve with fixed mass flux the 
depth of the solutions varies. 

Recently, in \cite{HDisp, HDisp13,DHSSG2016}, questions concerning the existence of wave 
wave solutions of fixed mean-depth with discontinuous vorticity, were addressed.  Physically, such 
studies are motivated by  the fact that wind blowing over a water surface induces a thin layer of 
high vorticity 
\cite{PhB74}. Mathematically, we note recent   deep analytical studies of flows with general 
discontinuous vorticity, see \cite{Cons_disp_12, CS11, Esc_bog_14, Bog_cal_14, BogAnc14}. 
This paper completes the analysis for the work \cite{HDisp, HDisp13,DHSSG2016} by proving the 
equivalence of the weak modified-height function formulation of the governing equations, to the 
weak Euler and stream function formulations. In \cite{DH1} the Euler, stream and modified-height 
formulations of the governing equations were proven to be equivalent in the sense of classical 
solutions. In the setting of the  classical Dubreil-Jacotin transformation and standard height 
function formulation,  this equivalence was recently proved by Constantin and Strauss  \cite{CS11}  
for solutions considered in a 
weak sense, for solutions in Sobolev spaces. More specifically, they 
consider weak solutions understood in the sense of distributions, and the solutions of weak wave 
formulation are considered in 
the Sobolev spaces. In \cite{MMAt_14} the equivalence was proved for strong solutions in Sobolev spaces that posses additional weak H\"older regularity. 
The result of equivalence for the standard height formulation was 
improved in \cite{VZ12}, where the authors consider weak solutions 
in H\"older spaces. The aim of this paper is to suitably adapt and generalise this result to prove 
the equivalence of the three weak formulations considering the modified-height function. 

In this paper, in section \ref{sec:gov:eq} we describe in detail the three  formulations: 
the  Euler equation; the  
stream formulation, and the modified-height formulation.  
In section \ref{sec:weak:form} we specify in detail the equations and 
boundary conditions of the weak form of the previous three formulations. 
And finally, in section \ref{sec:equiv:form}, we give the 
main result of equivalence between the three formulations. 
The equivalence between the weak 
stream and the weak modified-height formulation is proved for weak solutions with H\"older 
regularity where $0<\alpha\le 1$. However, since the equivalence between the weak velocity and the weak stream 
formulation is only proved  for weak solutions with H\"older regularity where $1/3<\alpha\le 1$, 
we can only affirm that the equivalence between the three formulations is satisfied for $1/3<\alpha
\le 1$. This equivalence of formulations is important, since in \cite{DHSSG2016}  
local bifurcation theory is applied to the modified-height formulation in order to prove the 
existence of small amplitude water waves on a fluid with discontinuous vorticity. 

\section{Standard Governing Equations}\label{sec:gov:eq}

We formulate the standard governing equations using Cartesian $(x,y)$-coordinates. 
These equations are defined in a frame which moves alongside the wave. Let $d>0$ 
be the depth of the undisturbed mass of water, and take $y=0$ to represent the level 
of the undisturbed water surface, then the flat bed is at $y=-d$. We assume the wave 
period is $2L$ and we denote by $\eta(x,t)$ the wave surface profile, which under physical 
assumptions satisfies that for any fixed time $t_0$ 
\[
	\int_{-L}^L\eta(x,t_0)dx=0.
\]
Without loss of generality, using scaling arguments, we work with $L=
\pi$ which will be more convenient. We are interested in travelling waves with a 
constant speed denoted by $c>0$ in the positive $x$-direction, then the velocity 
field takes the form $\left(u(x-ct,y),v(x-ct,y)\right)$ 
and the wave surface profile is given by $\eta(x-ct)$. One of the difficulties of 
this problem is that the wave profile $\eta$ is a free surface which 
is an unknown in the problem, then with the change of coordinates 
$(x-ct,y)\mapsto(x,y)$ we simplify the problem, obtaining now a time 
independent problem. We denote the fluid domain by $\overline{D_{\eta}}=\{(x,y)\in 
\mathbb R^2: -d\leq y\leq \eta (x)\}$, and the governing equations of the 
inviscid incompressible fluid is given by the mass conservation equation and Euler's 
equations together with the boundary conditions,
\begin{subequations}\label{Euler_subeqtions}
\begin{eqnarray}\label{masscon}
	u_x+v_y=0 &\mbox{ in } D_{\eta},\\
\label{euler'}
	(u-c)u_x+vu_y=-P_x &\mbox{ in } D_{\eta},\\ 
\label{euler'2}	(u-c)v_x+vv_y=-P_y-g &\mbox{ in } D_{\eta},\\
%
\label{ksc'}
	v=(u-c)\eta_x & \mathrm{on}\ y=\eta(x), \\
%
%
\label{bnddyn1}  
	P=P_{atm} & \mathrm{on}\ y=\eta(x),  \\
%
%
\label{kbc}
	v=0 & \mathrm{on}\ y=-d,
\end{eqnarray}
\end{subequations}
where $P=P(x,y)$ is the pressure, $g$ is the gravitational constant, and  $P_{atm}$ 
is the constant atmospheric pressure. The kinematic boundary 
condition \eqref{ksc'} express the fact that a particle on the free boundary 
remains there at all times; \eqref{bnddyn1} decouples the 
motion of the air from that of the water, and the last boundary condition \eqref{kbc}
assumes that the fluid does not penetrate the flat bed.  
The Eulerian governing equations for the gravity water wave problem are given by 
\eqref{Euler_subeqtions}, and for two-dimensional flows the vorticity is given by
\begin{equation}\label{vort}
	\omega=u_y-v_x.
\end{equation}
We assume also that the fluid does not contain stagnation points, that is %
\begin{equation}\label{umax} 
	u<c
\end{equation}
throughout the fluid. This means that the particles of the fluid move with less 
velocity than the wave speed, and physically, this assumption is valid for flows 
which are not near breaking, \cite{J90,Light}. 

The previous governing equations \eqref{Euler_subeqtions} can be reformulated in terms of the 
stream function $\psi$ which is directly related to $u, v$ by
\begin{equation}\label{psi1}
	\psi_y=u-c,\quad \psi_x=-v.
\end{equation}
This function is determined up to a constant. To fix the constant we consider 
$\psi=0$ on $y=\eta(x)$. We know from the boundary conditions  \eqref{ksc'} 
and \eqref{kbc} that $\psi$ is constant on both boundaries of $D_\eta$, and 
integrating \eqref{psi1} we obtain that $\psi=-p_0$ on $y=-d$, where 
\[
	p_0=\int_{-d}^{\eta(x)}(u(x,y)-c)dy
\] 
 is known as the relative mass flux, and thanks to \eqref{umax} we know that 
 for any given flow, $p_0$ is a fixed constant $p_0<0$. Integrating  \eqref{psi1}, we have that
\[
	\psi(x,y)=-p_0+\int_{-d}^{y}(u(x,s)-c)ds,
\]
and we can  see that $\psi$ is periodic in $x$, with period $2\pi$. 
Using \eqref{vort} and \eqref{psi1} we obtain that the stream function satisfies 
the equation
\[
	\Delta \psi=\omega,
\]
with 
$$\omega=\gamma(\psi/p_0),$$ 
where $\gamma$ is the vorticity function. Let 
\begin{equation*}
	\tilde\Gamma(p)=\int_{0}^{p}p_0\gamma(s)ds, \mbox{ for } -1\le p\le 0,
\end{equation*} 
then from Euler's equation \eqref{euler'} we obtain Bernoulli's law which states that 
\[
	E:= \frac{(u-c)^2+v^2}{2}+g(y+d)+P-\tilde\Gamma\left(\frac{\psi}{p_0}\right)
\]
is constant throughout the flow $\overline{D_\eta}$. We define  
$Q:=E-P_{atm}+gd$, then rewriting the governing equations in the moving frame 
in terms of the stream function, we obtain 
\begin{subequations}\label{Stream}
	\begin{eqnarray}\label{eqnsfirst}
		\Delta\psi = \omega  & \mathrm{in} & -d<y<\eta(x),\\
 		\lvert \nabla \psi \rvert ^2+2g(y+d)=Q  &\mathrm{on}& y=\eta(x), 
		\label{eulerorig}  
		\\ 
		\label{psi3}
		\psi=0 &\mathrm{on}& y=\eta(x),	
		\\ 
		\label{stokeslast}
		\psi=-p_0 &\mathrm{on}& y=-d.    
	\end{eqnarray}  
\end{subequations}
Furthermore, the condition which excludes 
stagnation points, \eqref{umax}, is equivalent to 
\begin{equation}\label{stag_psi}
	\psi_y<0.
\end{equation}
The main difficulties of solving the latter problem \eqref{Stream} are its nonlinear 
character and the fact that the free surface is unknown. 
To overcome this difficulty, we define the nonstandard semi-hodograph transformation,  
which was first introduced in \cite{DH1}, 
given by
\begin{equation}\label{hod}
q= x, \quad p= \frac{\psi (x,y)}{p_0}.
\end{equation}
This change of variables represents an isomorphism thanks to the assumption \eqref{umax}. The semi-hodograph 
transformation \eqref{hod} transforms the fluid domain $D_{\eta}$, with the 
unknown free boundary $\eta$, into the fixed semi-infinite rectangular domain 
$\overline R=\mathbb R\times [-1,0]$. 
%
We can now define the modified-height function in the $(q,p)$-variables,
\begin{equation}\label{hdef}
h(q,p)=\frac{y}d-p,
\end{equation}
where $y=y(q,p)$ is a function of the new $(q,p)$-variables. We assume 
that the modified-height function $h$ is even and $2\pi$-periodic on $q$, 
and by definition \eqref{hdef} and from \eqref{hod},  it satisfies
\begin{equation}\label{hint}
\int_{-\pi}^{\pi} h(q,0)dq=0.
\end{equation}
The 
modified-height function \eqref{hdef} was  introduced in \cite{ 
DH1}, where it was used to obtain existence results for rotational water waves 
of fixed mean-depth. This is a different approach from the approach taken in \cite{CS11}, 
where the authors fix the mass flux $p_0$ to prove the existence 
of solutions using local bifurcation. Here, as well as in in \cite{ HDisp, HDisp13, DH1} we fix 
the mean-depth.  
This is the more ideal physical approach, since it is easier to directly determine 
the mean-depth of a mass of water over a flat bed than the mass flux 
which is a more variable characteristic for any given flow. 
As was stated in \cite{DH1}, another difference 
in approaches from that of \cite{CS11} where the standard height function 
eliminates the depth $d$ from the problem, here  the 
modified-height function \eqref{hdef} allows us to introduce the depth 
$d$ into the problem. This is important, since in \cite{KoS1} was notice 
that fixing the mass flux does no fix the depth, since given a fix mass flux 
 there exists a bifurcation curve 
which has solutions with different depths. 
The semi-hodograph transformation \eqref{hod} 
transforms the stream function system of equations \eqref{eqnsfirst}-\eqref
{stokeslast} on an unknown domain with a free surface, into the following 
modified-height function system in a fixed domain
\begin{subequations}\label{H_0}
\begin{eqnarray}\label{hfirst}
\left(\frac{1}{d^2}+h^2_q\right) h_{pp}-2h_q(h_p+1)h_{pq}+(h_p+1)^2 h_{qq} 
+ \frac{\gamma(p)}{p_0}(h_p+1)^3 =0 & \mbox{ in } -1<p<0, \\ 
\label{hsecond}
\frac 1{d^2}+h_q^2+\frac{(h_p+1)^2}{p_0^2}[2gd(h+1)-Q]=0 &  \mbox
{ on } p=0,
\\ \label{hbnd}
h=0 &  \mbox{ on }  p=-1,
\end{eqnarray} 
\end{subequations}
where $h$ is even and $2\pi-$periodic in $q$ and \eqref{hint} holds. The condition 
which excludes stagnation points \eqref{stag_psi} is equivalent to 
\begin{equation}\label{hmax}
h_p+1>0,
\end{equation} 
and consequently the system \eqref{H_0} is a uniformly elliptic quasilinear partial 
differential equation with oblique nonlinear boundary condition.

\section{Weak formulations}\label{sec:weak:form}

In this section we describe the weak formulations which are associated to each 
of the formulations described in the previous section. 
These generalised formulations will give a meaning to solutions with weaker regularity 
than those of the formulations above, and the weak solutions 
we are interested in will be considered to be in H\"older spaces.

\subsection{Weak Euler equation} 
We can write the Euler equation 
\eqref{Euler_subeqtions} in the (weak) divergence form as 
\begin{subequations}\label{euler_weak}
\begin{eqnarray}
	\label{weul1}
	-cu_x+\left(u^2\right)_x+\left(uv\right)_y=-P_x\quad\qquad  \mbox{ in } D_{\eta},\\ 
	\label{weul2}
	-cv_x+\left(uv\right)_x+\left(v^2\right)_y=-P_y-g \quad\qquad\mbox{ in } D_{\eta},\\
	\label{weul3}
	u_x+v_y=0\quad\qquad \mbox{ in } D_{\eta},\\
	\label{weul4}
	v=0 \quad\, \mbox{on}\ y=-d,  \\
	\label{weul5} 
	v=(u-c)\eta_x \,\;\;\mbox{on}\ y=\eta(x),  \\
	\label{weul6}
	P=P_{atm} \,\;\;\mbox{on}\ y=\eta(x).
\end{eqnarray}
\end{subequations}  
In this formulation, the equations \eqref{weul1}--\eqref{weul3} will be understood 
in the sense of distributions, whereas the 
boundary conditions  \eqref{weul4}--\eqref{weul6} will be understood in the 
classical sense. The type of solutions of \eqref{euler_weak} we are 
interested in are solutions $u,v,P\in C^{0,\alpha}_{per}(\overline{D_{\eta}})$, 
where $\eta\in C^{1,\alpha}_{per}(\R)$, for some $\alpha\in (0,1]$. Here the 
\emph{per} subscript indicates that our solutions are even and $2\pi$-periodic 
in the $x$-variable. We assume also that the fluid has no stagnation 
points, so that $u<c$.

\subsection{Weak stream formulation}
Since the following  identity holds for regular enough functions $\psi$ 
and  $\gamma$, 
\begin{equation*}
	\big\{\psi_x\,\psi_y\big\}_x-{1\over 2}\big\{\psi_x^2-\psi_y^2\big\}_y-\big\{\tilde\Gamma(\psi/
	p_0)\big\}_y=\psi_y\left[\Delta\psi-\gamma(\psi/p_0)\right],
\end{equation*}
we can write the weak stream formulation 
as 
\begin{subequations}\label{stream_weak}
\begin{eqnarray}
	\label{wstr1}
	\big\{\psi_x\,\psi_y\big\}_x-{1\over 2}\big\{\psi_x^2-\psi_y^2\big\}_y-\big\{\tilde\Gamma
	(\psi/
	p_0)\big\}_y=0\quad\qquad  \mbox{ in } D_{\eta},\\ 
	\label{wstr2}
	\psi=-p_0 \quad\, \mbox{on}\ y=-d,  \\
	\label{wstr3} 
	\psi=0 \,\;\;\mbox{on}\ y=\eta(x),  \\
	\label{wstr4}
	\lvert \nabla \psi \rvert ^2+2g(y+d)=Q \,\;\;\mbox{on}\ y=\eta(x).
\end{eqnarray}
\end{subequations}  
Again this weak formulation \eqref{stream_weak} will give a meaning for 
solutions with weaker regularity. In this case we are interested in solutions 
$\psi\in C^{1,\alpha}_{per}(\overline{D_{\eta}})$, $\tilde\Gamma\in C^{0,\alpha}
([-1,0])$, and $\eta\in C^{1,\alpha}_{per}(\R)$, for some $\alpha\in 
(0,1]$, and the stream function has to satisfy the condition of there 
being no stagnation points, that is $\psi_y<0$.  The boundary conditions \eqref
{wstr2}--\eqref{wstr4}  are satisfied in the classical sense, and  the 
equation \eqref{wstr1} is satisfied in the sense of distributions. Notice that 
$\psi_x,\psi_y$ can be understood in the classical sense.

\subsection{Weak modified-height formulation}
We can rewrite the height equation in the divergence form 
\begin{subequations}\label{H}
\begin{align}\label{hfirst_0}
\left\{-\frac{1+d^2h^2_q}{2d^2(1+h_p)^2}+\frac{\Gamma(p)}{2d^2}\right\}_p+
\left\{\frac{h_q}{1+h_p} \right\}_q&=0  &\mbox{ in }  -1&<p<0, \\ \label{hsecond_0}
-\frac{1+d^2h^2_q}{2d^2(1+h_p)^2} -\frac{gd(h+1)}{p_0^2}+\frac{Q}
{2p_0^2}&=0 & \mbox{ on }  p&=0,
\\ \label{hbnd_0}
h&=0 &  \mbox{ on }  p&=-1.
\end{align} 
\end{subequations}
Here 
\Beqn
\Gamma(p)=2\int_0^p\frac{d^2\gamma(s)}{p_0}ds \quad\mbox{in } -1\leq p \leq 0.
\Eeqn
We understand by a solution of \eqref{H} a function  $h\in C^{1,\alpha}_{per}(\overline 
R)$, where $\Gamma\in C^{0,\alpha}_{per}([-1,0])$, for some $\alpha
\in(0,1]$. Here the \emph{per} subscript indicates that our solutions are 
even and $2\pi$-periodic in the $q$-variable. We assume also that the modified-height 
funtion satisfies $h_p+1>0$. The boundary conditions 
\eqref{hsecond_0} and 
\eqref{hbnd_0} are satisfied in the classical sense, and the equation \eqref
{hfirst_0} is satisfied in the sense of distributions. Thanks to the regularity 
considered for $h$, its derivatives $h_p,h_q$ are 
understood in the classical sense.  


\section{Equivalent formulation}\label{sec:equiv:form}
In this section we prove the main result which states the equivalence 
between the three weak formulations of the governing equations introduced in 
the previous section. In particular, 
we prove that the weak stream function system of equations, \eqref{stream_weak}, 
and the weak modified-height formulation, \eqref{H}, are equivalent with 
$\psi\in C^{1,\alpha}_{per}(\overline{D_{\eta}})$ and $h\in C^{1,\alpha}_{per}\left
(\overline{R}\right)$, for $0<\alpha \le 1$. On the other hand, the equivalence between 
the weak Euler equation \eqref{euler_weak},  and the weak stream function system of equations, 
\eqref{stream_weak}, are only proved for $1/3<\alpha\le  1$. This is why the result 
below is only proved for $\alpha$ between the latter values.

\begin{theorem}\label{equiv:form}
	Let $1/3<\alpha\le  1$. Then the following formulations of the governing equations 
	are  equivalent:
	\renewcommand{\theenumi}{\roman{enumi}}
	\begin{enumerate}
	\item[$(i)$] the weak Euler equation \eqref{euler_weak} with  
	\eqref{umax}, for  $\eta\in C^{1,\alpha}_{per}(\R)$, and 
	$u,v, P\in C^	{0,\alpha}_{per}(\overline{D_{\eta}})$;
	\item[$(ii)$] the weak stream  formulation \eqref{stream_weak} with 
	\eqref{stag_psi}, 
	for $\tilde{\Gamma}\in C^{0,\alpha}([-1,0])$,  $\eta
	\in C^{1,\alpha}_{per}(\R)$ and  $\psi\in C^{1,\alpha}_{per}(
	\overline{D_{\eta}})$;  
	\item[$(iii)$] the weak modified-height formulation \eqref{H} 
	together with \eqref{hmax}, for $\Gamma\in C^{0,\alpha}([-1,0])$ and  
	$h\in C^{1,\alpha}_{per}\left(\overline{R}\right)$.
	\end{enumerate}
\end{theorem}
%
%
\begin{proof}
	Let us prove first the equivalence between  the weak stream formulation $(ii)$ and 
	the weak modifed-height formulation $(iii)$ for $0<\alpha\le 1$.
	
	$\boxed{(ii)\Rightarrow (iii)}$ Let $\psi\in C^{1,\alpha}_{per}
	(\overline{D_{\eta}})$ satisfy \eqref{stream_weak} and \eqref{stag_psi}, 
	with $\tilde \Gamma\in C^{0,\alpha}([-1,0])$.  We recall the 
	semi-hodograph transformation given by 
	\begin{equation}\label{shod}
		(x,y)\mapsto (q,p)=\left(x, {\psi(x,y)\over p_0}\right), 
	\end{equation}
	is a bijection between $\overline{D_{\eta}}$ and $\overline{R}$ as a result of 
	\eqref{stag_psi}.  
	Let $h$ be the modified-height function $h(q,p)={y(q,p)\over d}-p$, then 
	the inverse mapping from $\overline{R}$ to $\overline{D_{\eta}}$ has 
	the form 
	\begin{equation}\label{inver_shod}
		(q,p)\mapsto (x,y)=\left(q, d[h(q,p)+p]\right).
	\end{equation}
	From the semi-hodograph transformation \eqref{shod} we obtain the 
	relations 
	\begin{equation}\label{partial:deriv}
		\partial_x=\partial_q-{h_q\over h_p+1}\partial_p, \quad
		\partial_y={1\over d(h_p+1)}\partial_p, \quad \mbox{and }\quad
		\partial_q=\partial_x-{\psi_x\over \psi_y}\partial_y,\quad
		\partial_p={p_0\over \psi_y}\partial_y,	
	\end{equation}
	and 
	\begin{equation}\label{der_h:u:v}
		\psi_x=-{p_0h_q\over h_p+1}, \quad 
		\psi_y={p_0\over d(h_p+1)},\quad \mbox{and }\quad
		h_q={-\psi_x\over d\psi_y},\quad 
		h_p={p_0\over d\psi_y}-1.
	\end{equation}
	The identities above should be regarded as a relation between  
	classical derivatives of a $C^1$-function with respect to the 
	$(x,y)$-variables and $(q,p)$-variables. Thanks to the relation between the 
	derivatives of $h$ and $\psi$, given by \eqref{der_h:u:v}, and since 
	$\psi\in C^{1,\alpha}_{per}(\overline{D_{\eta}})$ satisfies that  
	$\psi_y>0$ we have that $h\in C^{1,\alpha}_{per}
	(\overline{R})$. Now, let us prove that the boundary 
	conditions of the modified-height function 
	\eqref{hsecond_0}-\eqref{hbnd_0} hold. From \eqref{shod}, we 
	know that $p$ can be seen as 
	a function of the variables $x$ and $y$, 
	and  from \eqref{inver_shod}, $h$ can be regarded as a function of $q$ 
	and $p$. Moreover, the  
	It follows directly from \eqref{shod} and 
	\eqref{wstr2} that  
	\begin{equation*}
		p={\psi(x,y)\over p_0}=-1\;\mbox{ on } y=-d,
	\end{equation*}
	then 
	$$
		h(q,p)={y(q,p)\over d}+1=0\;\mbox{ on } p=-1.
	$$
	Thus, we have proved the boundary condition \eqref{hbnd_0}. On the 
	other hand,  thanks to \eqref{shod} and \eqref{wstr3} we have that 
	\begin{equation}\label{equiv_p_psi}
		p={\psi(x,y)\over p_0}=0\;\mbox{ on } y=\eta(x).
	\end{equation}
	From the relation \eqref{der_h:u:v} and thanks to \eqref{equiv_p_psi}, we 
	have that the boundary condition \eqref{hsecond_0} can be rewritten 
	as follows, and since we assume that the 
	boundary condition on the free surface \eqref{wstr4} for the stream 
	function is satisfied, we have that  
	\begin{equation*}
	\begin{array}{ll}
		\displaystyle	\left(\frac 1{d^2}+h_q^2\right)\frac{p^2_0}{(h_p+1)
		^2}+2gd(1+h+p)
		&
	 	\displaystyle=\frac{p_0^2h_q^2}{(h_p+1)
		^2}+\frac{p^2_0}{d^2(h_p+1)^2}+2gd\left(1+\frac y d\right)
		\smallskip\\
		&\displaystyle=\lvert \nabla 
		\psi \rvert ^2+2g(y+d)=Q \qquad  \mbox{ on } p=0.
	\end{array}	
	\end{equation*}
	Thus, \eqref{hsecond_0} is satisfied. Since the assumption of 
	there being no stagnation points for the stream function \eqref{stag_psi} 
	is satisfied and thanks to the relations \eqref{der_h:u:v}, we know that 
	$\psi_y={p_0\over d(h_p+1)}$, then the modified-height function 
	satisfies $h_p+1>0$, and \eqref{hmax} holds.
	
	Now, since we are considering weak solutions in $C^{1,\alpha}_{per}
	(\overline{R})$, to prove that the equation \eqref{hfirst_0} of the 
	modified-height function is satisfied, we prove it in the sense of 
	distributions. We need to prove that 
	\begin{equation}\label{weak_form_q_p}
	\int \int_R \left(-\frac{1+d^2h^2_q}{2d^2(1+h_p)^2}+\frac{\Gamma(p)}
	{2d^2}\right)\tilde\varphi_p+\left(\frac{h_q}{1+h_p} \right)\tilde
	\varphi_q\,dq\,dp=0\quad \mbox{for all }\;\tilde\varphi\in C^1_0(R).
	\end{equation}
	For any $\tilde\varphi$, let $\varphi\in C^1_0(D_{\eta})$ be given by 
	$\varphi(x,y)=\tilde\varphi\left(x,{\psi(x,y)/ p_0}\right)$ for all $(x,y)\in 
	D_{\eta}$. 
	Changing variables and from the relations \eqref{partial:deriv} and 
	\eqref{der_h:u:v} yields
	\begin{equation}
	\begin{array}{ll}\label{weak_form_cal}
	I 
	&=\displaystyle \int \int_{D_{\eta}} \left[\left(-{\psi_y^2\over 2 p_0^2}-
	{\psi_x^2\over 
	2p_0^2}+{\Gamma\left({\psi/ p_0}\right)\over 2d^2}\right){p_0\over 
	\psi_y}\varphi_y-{\psi_x\over p_0}\left(\varphi_x-{\psi_x\over \psi_y}
	\varphi_y \right)\right]{\psi_y\over p_0}dx\,dy  
	\smallskip\\
	&=\displaystyle \int\int_{D_{\eta}} \left[{1\over 2p_0^2}\left(-\psi_x^2-
	\psi_y^2\right)\varphi_y+ {\Gamma\left({\psi/ p_0}\right)\over 
	2d^2}\varphi_y-{\psi_x \psi_y\over p_0^2}\varphi_x+ {\psi_x^2\over p_0^2}\varphi_y\right] dx\,dy .
	\end{array}
	\end{equation}
	Multiplying \eqref{weak_form_cal} by $p_0^2$, and since $\tilde
	{\Gamma}(p)=\int_0^p p_0\gamma( s)ds={p_0^2\over 2d^2}\Gamma(p)
	$, 
	we have that
	\begin{equation}\label{weak_psi}
	p_0^2\,I =\displaystyle \int\int_{D_{\eta}}  \left[\tilde\Gamma\left({\psi/ p_0}\right)
	\varphi_y-
	\left(\psi_x \psi_y\right)\varphi_x+{1\over 2}\left
	(\psi_x^2-\psi_y^2\right)\varphi_y\right]dx\,dy.
	\end{equation}
	Since the stream function satisfies \eqref{wstr1} we have that $I=0$. Then 
	we have proved 	\eqref{weak_form_q_p}, and $(iii)$ holds.
	
	$\boxed{(iii)\Rightarrow (ii)}$ Let $h\in C^{1+\alpha}_{per}(\overline{R})
	$ satisfy \eqref{H} and \eqref{hmax}, with $\Gamma\in C^{0,\alpha}
	([-1,0])$. From the definition of $\psi$ and $\eta$,  we have that the 
	relation of the derivatives \eqref{partial:deriv} and \eqref{der_h:u:v}  are 
	still valid, and since $h\in C^{1+\alpha}_{per}(\overline{R})$, we have that  
	$\psi\in C^{1,\alpha}_{per}(\overline{D_
	{\eta}})$ and $\eta\in C^{1,\alpha}_{per}(\R)$. 
	Now, to recover $\psi$, we observe that the mapping $(q,p)\to (q, d(h(q,p)+p))$ 
	is a global bijection from $R$ onto $D_{\eta}$, because for $q$ fixed it is strictly 
	monotone and hence bijective. Moreover, the bijection is a global 
	$C^{1,\alpha}$-diffeomorphism. The inverse of this bijection is given by 
	$(x,y)\to (x,\psi(x,y)/p_0)$. 
	Now, let us prove that the boundary conditions 
	for the stream function  \eqref{wstr2}-\eqref{wstr4} 
	hold.
	 It follows directly from \eqref{shod} that  
	\begin{equation*}
		\psi(x,y)=-p_0\;\mbox{ on } y=-d,
	\end{equation*}
	and
	\begin{equation*}
		\psi(x,y)=0\;\mbox{ on } y=\eta(x),
	\end{equation*}
	and so $\psi$ satisfies the boundary conditions 
	\eqref{wstr2} and \eqref{wstr3}. From the definition of $h$, we have that  
\begin{equation}\label{yp}
y=d[h(x,p)+p].
\end{equation} 
	Differentiating  \eqref{yp} we get
	\begin{eqnarray}\label{wa1}
		y_x=0=d[h_q+h_pp_x+p_x] \Rightarrow p_x=-\frac{h_q}{1+h_p}, 
		\\ \label{wa2}
		y_y=1=d[h_pp_y+p_y] \Rightarrow p_y=\frac 1 {d(1+h_p)}.
	\end{eqnarray}
	Thanks to \eqref{der_h:u:v} and the previous relations \eqref{wa1} and \eqref{wa2},
	we can rewrite \eqref{wstr4} in terms of the modified-height function as 
	follows, and since \eqref{hsecond_0} is satisfied, we have that
 	\begin{equation*}
	\begin{array}{ll}
	 	\displaystyle\lvert \nabla \psi \rvert ^2+2g(y+d)&\displaystyle=\frac{p_0^2h_q^2}{(h_p+1)
		^2}+\frac{p^2_0}{d^2(h_p+1)^2}+2gd(1+\frac y d)
		\smallskip\\ 
		&\displaystyle=\left(\frac 1{d^2}+h_q^2\right)\frac{p^2_0}{(h_p+1)^2}+2gd(1+h+p)
		=Q  \quad \mbox{ on } y=\eta(x).
	\end{array}
	\end{equation*}	
	Thus, the boundary condition \eqref{wstr4} follows.  
	On the other hand, thanks to the relation \eqref{der_h:u:v} we have that $\psi_y={p_0\over 
	d(h_p+1)}$, and since $h$ satisfies 
	\eqref{hmax} and $p_0<0$, we have that the assumption of there being 
	no stagnations 	points for the stream function, \eqref{stag_psi} holds.  
	Now, let us prove that the equation \eqref{wstr1} is 
	satisfied in the distribution sense. To do this we have to prove that   
	 \begin{equation}\label{I_zero}
	 	\displaystyle \int\int_{D_{\eta}}  \tilde\Gamma\left({\psi/ p_0}\right)
		\varphi_y-\left(\psi_x \psi_y\right)\varphi_x+{1\over 2}\left(
		\psi_x^2-\psi_y^2\right)\varphi_ydxdy=0\quad \mbox{ for all }\;\varphi\in C^1_0(D_{\eta}).
	 \end{equation}
	For any $\varphi\in C^1_0(D_{\eta})$, let $\tilde\varphi\in C^1_0(R)$ be 
	given by $\tilde\varphi(q,p)=\varphi\left(q,d[h(q,p)+p]\right)$ for all $(q,p)
	\in R$. By changing variables in the integral in \eqref{I_zero} which is equal to 
	\eqref{weak_psi}, and following the arguments above from the bottom up, we can rewrite \eqref{weak_psi} 
	as \eqref{weak_form_q_p}. But \eqref{weak_form_q_p} is valid, as a 
	consequence of  \eqref{hfirst_0}. Hence, we have proved that $(ii)$ 
	holds.
	
	Although the details of the proof of the equivalence between $(i)$ and $(ii)$ follows as in 
	\cite{VZ12}, since the precise composition of the modified-height function plays no role 
	in the equivalence considerations, for the sake of completeness we present an outline of the 
	proof considering the stream formulation presented in this paper -- full details may be found 
	in \cite{VZ12}. 
	
%
	
	$\boxed{(i)\Rightarrow (ii)}$ 
	Since $(i)$ holds, then 
	 $u,v\in C^{0,\alpha}_{per}
	(\overline{D}_{\eta})$ and $\eta\in C^{1,\alpha}_{per}(\R)$. From the definition of 
	$\psi$, \eqref{psi1}, we have that $\psi\in C^{1,\alpha}_{per}(\overline{D_
	{\eta}})$, which is unique up to a constant. It is clear that \eqref{weul4} 
	and 	\eqref{weul5} imply \eqref{wstr2} and \eqref{wstr3}. Using the 
	definition of $\psi$ we rewrite \eqref{weul1} and \eqref{weul2} in the 
	weak distributional form (with the first derivatives understood in the 
	classical sense),
	\begin{subequations}\label{weak_psi_proof}
	\begin{align}\label{weak_psi_proof_1}
		\left(\psi_y^2\right)_x-\left(\psi_x\psi_y\right)_y=-P_x\;\mbox{ in } D_
		{\eta},\\
		\label{weak_psi_proof_2}
		-\left(\psi_x\psi_y\right)_x+\left(\psi_x^2\right)_y=-P_y-g\;\mbox{ in } D_
		{\eta}.
	\end{align}
	\end{subequations}
	Let us define
	\begin{equation}\label{def_F_proof}
		F:=P+{1\over 2}|\nabla \psi|^2+gy\;\mbox{ in } D_{\eta},
	\end{equation}
	then it follows from \eqref{weak_psi_proof_1} that  the derivatives of $F$ 
	in sense of distributions are given by  
	\begin{eqnarray}\label{F_x_00}
		F_x={1\over 2}\left(\psi_x^2-\psi_y^2\right)_x +\left(\psi_x\psi_y\right)_y,
		\\ \label{F_y_00}
		F_y=\left(\psi_x\psi_y\right)_x-{1\over 2}\left(\psi_x^2-\psi_y^2\right)_y.
	\end{eqnarray}
	Now we prove that the equation \eqref{wstr1} holds, which is \eqref{F_y_00} in the sense of 
	distributions if we show that there exists a function 
	$\tilde\Gamma\in C^{0,\alpha}\left([p_0,0]\right)$ such that 
	\begin{equation}\label{concl_proof_i_ii_0}
		F(x,y)=\tilde\Gamma(\psi(x,y)/p_0) \;\mbox{ for all } (x,y)\in D_{\eta},
	\end{equation}
	where $\tilde\Gamma({p})=\int_{0}^{p}p_0\gamma(s)ds$. 
	Let $\tilde F: R\to \R$ be given 
	by $\tilde{F}({q},{p})=F\left({q},y({q},{p})\right)$ in $R$, which is 
	equivalent to $F(x,y)=\tilde F\left(x,\psi(x,y)/p_0\right)$ in $\overline{D}_{\eta}$. 
	Then \eqref{concl_proof_i_ii_0} is equivalent to 
	\begin{equation}\label{concl_proof_i_ii}
		\tilde F({q},{p})=\tilde\Gamma({p}),
	\end{equation}	
	for some $\tilde\Gamma\in C^{0,\alpha}([-1,0])$. To prove \eqref
	{concl_proof_i_ii}, we have to see if 
	\begin{equation}\label{concl_proof_i_ii_def}
		\int\int_{R}\tilde F\tilde{\varphi}_{q} d{q}d{p}=0 \;\mbox{for all } \tilde\varphi\in C^1_0(R).
	\end{equation}	
	For any $\tilde{\varphi}$, let $\varphi\in C_0^1(D_{\eta})$ be given by 
	$\varphi(x,y)=\tilde{\varphi}(x,\psi(x,y)/p_0)$ for all $(x,y)\in D_{\eta}$. 
	Changing variables in the integral \eqref{concl_proof_i_ii_def} we obtain 
	\begin{equation}\label{concl_proof_i_ii_def_2}
		\int\int_{D_{\eta}}F\left(\psi_y\varphi_x-\psi_x\varphi_y\right)dxdy=0.
	\end{equation}	
	Our aim is to prove \eqref{concl_proof_i_ii_def_2} for all 
	$\varphi\in C_0^1(D_{\eta})$. We define $V:=D_{\eta}$, and for  
	$\varphi \in C^1_0(D_{\eta})$ arbitrary, we denote $K:=\supp \varphi$. 
	Let $\eps_0:=dist (K, \R^2\setminus V)/2$, then for $0<\eps<\eps_0$, 
	we denote $V_{\eps}:=\left\{x\in V: dist (x,\R^2\setminus D_{\eta})>
	\eps\right\}$. Let $\rho^{\eps}$ be a mollifier defined in $V^{\eps}$, and let   
	$F^\eps:=F\ast \rho^{\eps}$ be defined in $V^{\eps}$. 
	Then we can write \eqref{concl_proof_i_ii_def_2} as 
	$$		
		 \begin{array}{ll}
		\displaystyle\int \int_{D_{\eta}}F& \!\!\!\!\!\big(\psi_y\varphi_x \displaystyle-\psi_x
		\varphi_y\big)dxdy=\\
		&\displaystyle\int\int_{K}\left(F\psi_y- F^\eps\psi^\eps_y\right)\varphi_x-\left(F\psi_x
		-F^\eps\psi^\eps_x\right)\varphi_ydxdy
		+\int\int_{K}F^\eps\psi_y^{\eps}\varphi_x-F^\eps\psi_x^{\eps}\varphi_ydxdy.
		\end{array}
	$$
	Thanks to \cite[Lemma 4.2]{VZ12}, we have some estimates of the norm 
	of $F^{\eps}$  and $\psi^{\eps}$ given in terms of the norm of $F$ and $\psi$. 
	Thanks to these estimates we obtain that  
	$$
		\int\int_{D_{\eta}}F\left(\psi_y\varphi_x-\psi_x\varphi_y\right)dxdy\le 
		C_1\eps^{\alpha}+C_2\eps^{2\alpha}+C_3\eps^{3\alpha-1},
	$$
	then if $\alpha>1/3$ and taking limits as $\eps$ goes to zero, we obtain 
	\eqref{concl_proof_i_ii_def_2}, and $(ii)$ holds.	
	
	$\boxed{(ii)\Rightarrow (i)}$ Let $u, v$ be defined by \eqref{psi1} up 
	to a constant, and the pressure by 
	$$
		P:=-{1\over 2}|\nabla \psi|^2-gy+\tilde\Gamma(\psi/p_0)\mbox{ in } 
		\overline{D}_{\eta}. 
	$$
	Then, $u,v,P\in C^{0,\alpha}_{per}(\overline{D}_{\eta})$. The definition of 
	$u$ and $v$ implies \eqref{weul3}, and \eqref{wstr2} and \eqref{wstr3}
	imply  \eqref{weul4} and \eqref{weul5}. Using the definition of $u,v$ and 
	$P$ mentioned above, we can rewrite \eqref{weul1} and \eqref{weul2} as 
	\begin{subequations}\label{rew_eul1_eul2}
	\begin{align}\label{rew_eul1_eul2_a}
	\tilde\Gamma(\psi/p_0)_x={1\over 2}\left(\psi_x^2-\psi_y^2\right)_x+\left
	(\psi_x\psi_y\right)_y\mbox{ in } D_{\eta},\\
	\label{rew_eul1_eul2_b}
	\tilde\Gamma(\psi/p_0)_y=\left(\psi_x\psi_y\right)_x-{1\over 2}\left
	(\psi_x^2-\psi_y^2\right)_y \mbox{ in } D_{\eta}.
	\end{align}
	\end{subequations}
	However, \eqref{rew_eul1_eul2_b} is exactly \eqref{wstr1}, which we 
	are assuming to hold. We just need to prove that \eqref{rew_eul1_eul2_a} holds, 
	and to do this we will prove that \eqref{rew_eul1_eul2_a} is a 
	consequence of \eqref{rew_eul1_eul2_b}. Let us define $F=
	\tilde{\Gamma}(\psi/p_0)$, then \eqref{rew_eul1_eul2_a} is equivalent to proving that   
	\begin{equation}\label{concl_proof_i_ii_def_3}
		\int\int_{D_{\eta}} F\varphi_x-{1\over 2}\left(\psi_x^2-\psi_y^2\right)\varphi_x
		-(\psi_x\psi_y)\varphi_ydxdy=0 \quad \mbox{for all } \varphi\in C^1_0(D_{\eta}).
	\end{equation}
	We define $V:=D_{\eta}$, and for  
	$\varphi \in C^1_0(D_{\eta})$ arbitrary, let $K:=\supp \varphi$. We 
	define $\eps_0:=dist (K, \R^2\setminus V)/2$, then for $0<\eps<\eps_0$, 
	let $V_{\eps}:=\left\{x\in V: dist (x,\R^2\setminus D_{\eta})>
	\eps\right\}$, and let $\rho^{\eps}$ be a mollifier defined in $V^{\eps}$, then we 
	consider $F^\eps:=F\ast \rho^{\eps}$ defined in $V^{\eps}$. 
	We can rewrite \eqref{concl_proof_i_ii_def_3} as 
	\begin{equation}\label{concl_proof_i_ii_def_4}
	\begin{array}{lll}
		\displaystyle \int\int_{D_{\eta}} &\!\!\!\!\!\!F&\displaystyle\!\!\!\!\!\!\!\!\!\!\!\varphi_x-{1\over 2}\left(\psi_x^2-\psi_y^2\right)\varphi_x
		-(\psi_x\psi_y)\varphi_ydxdy
		\\\smallskip
		&=&\displaystyle\int\int_{K} \left[F-F^{\eps}\right]\varphi_x-\left[{1\over 2}
		(\psi_x^2-\psi_y^2)-{1\over 2}\left(\left(\psi_x^{\eps}\right)^2\left(\psi_y^{\eps}\right)^2\right)
		\right]\varphi_x
		\\\smallskip 
		&&\displaystyle-\left[(\psi_x\psi_y) -\left(\psi_x^{\eps}\psi_y^{\eps}\right)\right]\varphi_ydxdy
		\\\smallskip 
		&&\displaystyle+\int\int_{K} F^{\eps}\varphi_x-{1\over 2}\left(\left(\psi_x^{\eps}\right)^2\left(
		\psi_y^{\eps}\right)^2\right)\varphi_x-\left(\psi_x^{\eps}\psi_y^{\eps}\right)\varphi_ydxdy.
	\end{array}
	\end{equation}
	Again, thanks to \cite[Lemma 4.2]{VZ12}, we have some estimates of the norm 
	of $F^{\eps}$  and $\psi^{\eps}$ given in terms of the norm of $F$ and $\psi$. 
	Thanks to these estimates we obtain that  
	$$
		\int\int_{D_{\eta}} F\varphi_x-{1\over 2}\left(\psi_x^2-\psi_y^2\right)\varphi_x
		-(\psi_x\psi_y)\varphi_ydxdy\le 
		C_1\eps^{\alpha}+C_2\eps^{2\alpha}+C_3\eps^{3\alpha-1}.
	$$
	Since $\alpha>1/3$, then taking limits as $\eps$ goes to zero, we obtain 
	\eqref{concl_proof_i_ii_def_3}. Thus, we have proved that $(i)$ holds.
\end{proof}

\end{document}